\documentclass[a4paper,11pt]{article}
\usepackage{a4wide}
\usepackage[latin1]{inputenc}
\usepackage{amsfonts}
\usepackage{verbatim}
\usepackage{mathrsfs}
\usepackage{amsmath}
\usepackage{amssymb}
\usepackage{amsthm}
\usepackage[T1]{fontenc}
\usepackage{graphicx}
\newcommand{\ud}{{ \, \mathrm{d}}}
\newcommand{\nd}{\stackrel{\mathrm{def}}{=}}
\newcommand{\lla}{\left\langle}
\newcommand{\rra}{\right\rangle}

\newtheorem{Theorem}{Theorem}[section]
\newtheorem{Definition}[Theorem]{Definition}
\newtheorem{Proposition}[Theorem]{Proposition}
\newtheorem{Lemma}[Theorem]{Lemma}

\newtheorem{Notation}[Theorem]{Notation}

\author{G. Fabbri\footnote {School of Mathematics and Statistics, UNSW, Sydney, e-mail: fabbri@maths.unsw.edu.au, supported by the ARC Discovery project DP0558539.},
 \; B. Goldys\footnote {School of Mathematics and Statistics, UNSW, Sydney, e-mail: goldys@maths.unsw.edu.au, supported by the ARC Discovery project DP0558539.}}
\title{An LQ problem for the heat equation on the halfline with Dirichlet boundary control and noise}
\begin{document}
\maketitle

\begin{abstract}
\bigskip

We study a linear quadratic problem for a system governed by the heat equation on a halfline with Dirichlet boundary control and Dirichlet
boundary noise. We show that this problem can be reformulated as a stochastic evolution equation in a
certain weighted $L^2$ space. An appropriate choice of weight allows us to prove a stronger regularity for the boundary
terms appearing in the infinite dimensional state equation. The direct solution of the Riccati equation related to the
associated non-stochastic problem is used to find the solution of the problem in feedback form and to write the value
function of the problem.

\bigskip



\end{abstract}

\section{Introduction}
In this paper we are concerned with a linear quadratic control problem for a heat equation on the halfline $[0,\infty)$
with Dirichlet boundary control and boundary noise. More precisely, for fixed $0\leq \tau <T$, we deal with the
equation
\begin{equation}
\label{eq:pde-state-equation}
\left \{
\begin{array}{ll}
\frac{\partial}{\partial t} y(t,\xi) = \frac{\partial^2}{\partial \xi^2} y(t,\xi) & t\in [\tau,T], \; \xi>0,\\
y(t,0)= u(t) + \dot{W}(t) &  t\in [\tau,T], \\
y(\tau,\xi)= x_0(\xi) & \xi>0.
\end{array}
\right .
\end{equation}
where $W$ is a one dimensional Brownian motion and $u$ is a square-integrable control. Let us recall that a
deterministic boundary control problem
\begin{equation}
\label{eq:pde-state-equation2} \left \{
\begin{array}{ll}
\frac{\partial}{\partial t} z(t,\xi) = \frac{\partial^2}{\partial \xi^2} z(t,\xi) & t\in [\tau,T], \; \xi>0,\\
z(t,0)= u(t) &  t\in [\tau,T], \\
z(\tau,\xi)= x_0(\xi) & \xi>0,
\end{array}
\right .
\end{equation}
is well understood, see for example \cite{BDDM2}, \cite{Lasiecka80}). Denoting by $A_0$ the Dirichlet Laplacian in
$L^2(0,\infty)$ and by $D$ the Dirichlet map (defined as $D_{\lambda_0}$ in (\ref{dir}) below), we can rewrite (\ref{eq:pde-state-equation2}) in the form
\[z(t)=e^{tA_0}x_0+\left(\lambda_0-A_0\right)\int_0^te^{(t-s)A_0}Du(s)ds\, ,\]
and it is easy to show that $z(t)\in L^2(0,\infty)$ for all $t\ge 0$. Therefore, the process
\begin{equation}\label{e3}
X(t)=e^{tA_0}x_0+\left(\lambda_0-A_0\right)\int_0^te^{(t-s)A_0}Du(s)ds+\left(\lambda_0-A_0\right)\int_0^te^{(t-s)A_0}DdW(s)
\end{equation}
 seems to be a good candidate for a solution
to (\ref{eq:pde-state-equation}). However, it was shown in \cite{DaPratoZabczyk93} that the process $X$ is not
$L^2$-valued. More precisely, it was shown that the solution to (\ref{eq:pde-state-equation})  considered on a finite
interval and for $u=0$, when rewritten in the form (\ref{e3}), is well defined in a negative Sobolev space
$H^{-\alpha}$ for $\alpha>\frac{1}{2}$ only. It is easy to see that the same conclusion holds in the case of halfline. Then it was shown in \cite{AlosBonaccorsi02}, see also
\cite{BonaccorsiGuatteri02}, that the process $X$ can be defined pointwise on $(0,\infty)$ and it takes values in a weighted space $L^2\left(0,\infty
;\xi^{1+\theta}d\xi\right)$. This fact was used to study some properties of the process $X$ (in fact in the
aforementioned papers more general nonlinear equations are studied) but the problem is not reformulated as a stochastic
evolution equation in $L^2\left(0,\infty ;\xi^{1+\theta}d\xi\right)$ and therefore advantages of using the weighted space are
somewhat limited.
\par\noindent
Following the idea of Krylov \cite{Krylov01} we introduce the weighted spaces $\mathcal H_\rho=L^2\left((0,\infty);\rho(\xi)d\xi\right)$, where for $\theta\in (0,1)$ we have
\[\rho(\xi)=\xi^{1+\theta}\quad\mathrm{or}\quad\rho(\xi)=\min\left(1,\xi^{1+\theta}\right),\quad \xi\ge 0.\]
It was proved in \cite{Krylov99} and \cite{Krylov01} that the
Dirichlet Laplacian $A_0$ defined on $L^2(0,\infty)$ extends to a generator $A$ of an analytic semigroup $\left(e^{tA}\right)$ on $\mathcal H_\rho$.
We will show that the Dirichlet map takes values in $\mathrm{dom}\left((-A)^\alpha\right)$
for a certain $\alpha>\frac{1}{2}$ and therefore equation (\ref{e3}), when considered in $\mathcal H_\rho$, can be given a form
\[X(t)=e^{tA}x_0+\int_0^te^{(t-s)A}\left(\lambda_0-A\right)Du(s)ds+\int_0^te^{(t-s)A}\left(\lambda_0-A_0\right)DdW(s)\]
that is, we will study a controlled evolution equation
\begin{equation}
\label{eq:state-equation-intro} \left \{
\begin{array}{l}
\ud x(t) = \left (  A x(t) + Bu(t) \right ) \ud t + B\ud W(t)\\
x(\tau)=x_0\in\mathcal{H}_\rho
\end{array}
\right .
\end{equation}
for $B=(\lambda_0 - A)D$.
 This fact is a starting point for our analysis of the linear
quadratic control problem (\ref{eq:pde-state-equation}). We will demonstrate that the control problem (\ref{eq:state-equation-intro}) when considered in the space $\mathcal H_\rho$ can be solved using classical by now techniques presented, for example, in \cite{BDDM2}. Let us emphasize that while focus of this paper is on the most interesting case of boundary control and boundary noise a more general control problem
\begin{equation}
\label{eq:state-equation-intro_g} \left \{
\begin{array}{l}
\ud x(t) = \left (  A x(t) + Bu_1(t)+v(t)\right ) \ud t + B\ud W(t)+dW_1(t)\\
x(\tau)=x_0\in\mathcal{H}_\rho
\end{array}
\right .
\end{equation}
with spatially distributed noise $W_1$ and control $v$ might be easily considered using the same technique.
\par\noindent
Let us note that if the boundary conditions are of Neumann type then the analogue of equation
(\ref{eq:pde-state-equation}) has a solution in $L^2(0,\infty)$ and has been studied intensely (also for more general
parabolic equations with boundary noise), see for example \cite{DaPratoZabczyk93}, \cite{Maslowski95}, \cite{DebusscheFuhrmanTessitore07}, \cite{dunmas}.
\par\medskip\noindent
We study the linear quadratic problem characterized by the cost functional
\begin{equation}
\label{eq:functional-intro}
J(\tau, x_0, u)=\mathbb{E} \left [ \int_\tau^T |Cx(t)|_Y^2 + |u(t)|^2_\mathbb{R} \ud t + \lla G x(T) ,x(T) \rra_{\mathcal{H}_\rho} \right ]
\end{equation}
and governed by a state equation of the form (\ref{eq:state-equation-intro}). the operator $C$ that appears
in (\ref{eq:functional-intro}) is in  $\mathcal{L}(\mathcal{H_\rho}; Y)$ for a certain Hilbert space $Y$ and
$G\in\mathcal{L}(\mathcal{H_\rho};\mathcal{H_\rho})$ is symmetric and positive. The direct solution of the Riccati equation
related to a linear quadratic problem driven by a stochastic equation different from ours was studied in the Neumann case
(non-weighted setting) in \cite{Flandoli86} (see also \cite{Ahmed81} and \cite{DaPrato84} for the control inside the
domain case ($\alpha=1$)). Our approach is different from the one used in the aforementioned works since we directly use the solution
of the Riccati equation for the ``associated'' deterministic problem.
\par\noindent
The deterministic linear quadratic problem \textit{associated} to ours is that characterized by the state equation
\[\dot{x}(t) = A x(t) + Bu(t)\]
 and the functional
\begin{equation}\label{det0}
\int_\tau^T \left(|Cx|_Y^2 + |u|^2_\mathbb{R}\right) \ud t + \lla G x(T) ,x(T)\rra_{\mathcal H_\rho}.
\end{equation}
It is well known, see \cite{BDDM2} and Section 3 below, that the solution to the linear quadratic problem given above is determined by the operator-valued function $P:[0,T]\to\mathcal L\left(\mathcal H_\rho, \mathcal H_\rho\right)$ which solves the so-called Riccati equation
\begin{equation}
\label{eq:riccati-intro}
\left \{
\begin{array}{l}
P'(t) = -A^*P(t) - P(t)A^* - C^*C + P(t)ABB^*A^*P(t)\\
P(T)= G.
\end{array}
\right .
\end{equation}
Such a problem has been intensely studied (see \cite{BDDM2} and \cite{LasieckaTriggiani00} and the references therein). We will refer in particular to the direct solution approach and we will use the formalism introduced in Section 2.2. of \cite{BDDM2}. We show that the Riccati equation (\ref{eq:riccati-intro}) has a unique solution $P$ in the space $C_{s,\alpha}([0,T];\Sigma(\mathcal{H}))$ (see Definition \ref{def:Csalphabeta}). Let us note that in the deterministic case the minimum of the cost functional \eqref{det0} is given by $\left \langle P(\tau) x_0, x_0 \right \rangle_{\mathcal{H}}$.
\par\noindent
In the study of the problem with boundary noise some of the tools and the results of the deterministic case, as the properties of the elements of $C_{s,\alpha}\left([0,T];\Sigma\left(\mathcal{H_\rho}\right)\right)$ and the solution of (\ref{eq:riccati-intro}), are still useful. It is possible to express the value function and the optimal feedback in terms of $P$. A term due to the noise appears in the expression of the minimal cost and we have that (Theorem \ref{th:optimalfeedback}):
\begin{multline}
V(\tau,x_0) = \inf_{u \in \mathcal{U}_\tau} J(\tau, x_0, u) =\\
=\lla P(\tau) x_0 ,x_0 \rra + \int_\tau^T \frac{1}{2} \left \langle ((\lambda_0 -A) D(1)), P(s)((\lambda_0 -A) D(1)) \right\rangle_{\mathcal{H}} \ud s.
\end{multline}
\section{The heat equation in $\mathcal{H_\rho}$}
\label{sec:preliminary}
\subsection{Notation}
We will work in a weighted space $\mathcal H_\rho=L^2([0,\infty);\, \rho(\xi)d\xi)$, where either $\rho(\xi) = \xi^{1+\theta}\wedge 1$ or $\rho(\xi) = \xi^{1+\theta}$ for some $\theta \in (0,1)$ and $\xi\geq 0$. All the results proved in the sequel are valid for both weights and therefore, in order to simplify notations \textbf{we will use the same notation $\mathcal{H}=\mathcal H_\rho$ for both weights}. Let us recall that $f\in\mathcal{H}$ if and only if
\[\int_0^{\infty} f^2(\xi) \rho(\xi) \ud \xi <\infty\]
and $\mathcal H$ is a Hilbert space with the scalar product
\[\left \langle \phi , \psi \right\rangle_{\mathcal{H}} = \int_0^{\infty} \phi(\xi) \psi(\xi) \rho(\xi) \ud \xi \qquad \text{for all } \phi, \psi \in \mathcal{H}.\]
Given $\lambda>0$, the Dirichlet map $D_{\lambda}$ is defined as follows:
\begin{equation}\label{dir}
D_{\lambda} (a) = \phi \; \Longleftrightarrow \; \left \{
\begin{array}{l}
(\lambda - \partial_x^2) \phi(\xi)=0 \qquad \text{for all } \xi>0\\
\phi(0)=a
\end{array}
\right .
\end{equation}
so $D_{\lambda} (a) = a \psi_{\lambda}$ where
\begin{equation}
\label{eq:formadilambda}
\left \{
\begin{array}{l}
\psi_\lambda \colon \mathbb{R}^+ \to \mathbb{R}\\
\psi_\lambda \colon \xi \mapsto e^{-\sqrt{\lambda}\xi}
\end{array}
\right .
\end{equation}
Clearly $\psi_\lambda\in\mathcal H$.\\
It is well known that for every $x_0\in L^2(0,\infty)$ the solution $y$ to the heat equation with zero Dirichlet boundary condition
\[
\left \{
\begin{array}{ll}
\frac{\partial}{\partial t} y(t,\xi) = \frac{\partial^2}{\partial \xi^2} y(t,\xi) & t>0, \; \xi>0,\\
y(t,0)= 0 &  t\ge 0, \\
y(0,\xi)= x_0(\xi) & \xi>0.
\end{array}
\right.\]
 is given by the following well known expression
\begin{equation}
\label{eq:solution-convolution} y(t,\xi)=\int_0^\infty k(t,\xi,\eta) x_0(\eta) \ud \eta
\end{equation}
where
\begin{equation}
\label{eq:kernel}
k(t,\xi,\eta) = \frac{1}{\sqrt{4\pi t}} \left ( e^{-\frac{(\xi-\eta)^2}{4t}} - e^{-\frac{(\xi+\eta)^2}{4t}}
\right ), \;\; \eta,\xi\ge 0.
\end{equation}
This formula defines the corresponding heat semigroup $T(t)x_0=y(t)$ in $L^2(0,\infty)$.
It is also well known that $(T(t))$ is a symmetric $C_0$-semigroup of contractions on $L^2(0,\infty)$.
\subsection{Properties of the heat semigroup on $\mathcal{H}$}
\begin{Proposition}
\label{pr:propertiessemigroup} For each of the weights $\rho(\xi)$ considered above, the heat semigroup $(T(t))$ extends
to a bounded $C_0$ semigroup $\left(e^{tA}\right)_{t\geq 0}$ on $\mathcal{H}$ with generator $A\colon \mathrm{dom}(A) \to
\mathcal{H}$. The semigroup $\left(e^{tA}\right)_{t\geq 0}$ is analytic.
\end{Proposition}
\begin{proof}
\emph{The case $\rho(\xi)=\xi^{1+\theta}$: $\mathcal{H}= L^2\left([0,\infty),\xi^{1+\theta}d\xi\right)$}.\\
Let $f\in L^2(0,\infty)$. Then by Theorem 2.5 in \cite{Krylov01} there exists $C>0$ independent of $f$ and such that
\[\left|e^{tA}f\right|_{\mathcal H}\le C|f|_{\mathcal H},\quad t\ge 0.\]
Since $L^2(0,\infty)$ is dense in $\mathcal H$, $e^{tA}$ can be extended to $\mathcal H$ and the strong continuity
follows by standard arguments. Let $A_0$ be the generator of $(T(t))$ in $L^2(0,\infty)$ and let $\mathcal D=\mathrm{dom}\left(A_0\right)\cap\mathcal H\subset\mathcal H$. Clearly
\[e^{tA}\mathcal D\subset\mathcal D,\quad t\ge 0\]
and $\mathcal D$ is dense in $\mathcal H$. Therefore $\mathcal D$ is a core for the generator $A$ of $\left(e^{tA}\right)$ in $\mathcal H$. If $f\in\mathcal D$
then
\[Ae^{tA}f=\frac{\partial^2}{\partial\xi^2}T(t)f\]
and again by Theorem 2.5 in \cite{Krylov01} we have
\[\left|\frac{\partial^2}{\partial\xi^2}T(t)f\right|_{\mathcal H}\le\frac{C}{t}|f|_{\mathcal H}.\]
Since $\mathcal D$ is a core for the generator $A$ in $\mathcal H$, the above estimate can be extended to any $f\in\mathcal H$ and therefore
\[\left|Ae^{tA}f\right|_{\mathcal H}\le \frac{C}{t}|f|_{\mathcal H}\quad f\in\mathcal H.\]
The last inequality is equivalent to the analyticity of the semigroup $\left(e^{ta}\right)$ in $\mathcal H$.
follows.\\
\emph{The case $\rho(\xi)=1 \wedge\xi^{1+\theta}$: $\mathcal{H}=L^2\left([0,\infty),1\wedge \xi^{1+\theta}d\xi\right)$}.\\
Let $x\in C_0^\infty(0,\infty)$ and $t\le T$. Then the functions $x_1=xI_{[0,1]}$ and
$x_2$ are in $L^2(0,\infty)$ and $\mathcal H_\rho$ for both weights $\rho$. It follows that
\begin{equation}
\label{eq:C0secondacase}
\begin{aligned}
|(T(t)x |_{\mathcal{H}} &\leq |T(t) (\chi_{[0,1]}x) |_{\mathcal H} +
|T(t) (\chi_{(1,+\infty)}x) |_{\mathcal H}\\
&\leq |T(t) (\chi_{[0,1]}x) |_{L^2_{\xi^{1+\theta}}} + |T(t) (\chi_{(1,+\infty)}x) |_{L^2(0,+\infty)}\\
&\leq C|x|_{\mathcal H}
\end{aligned}
\end{equation}
for a certain $C>0$. The fact that $C$ does not depend on $t\le T$ is a consequence of the $C_0$ property of $T_t$ on
$L^2_{\xi^{1+\theta}}$ (showed in the first part of the proof) and on $L^2(0,\infty)$. Therefore $(T(t))$ has an
extension to a semigroup $\left(e^{tA}\right)$ on $\mathcal H$ and the $C_0$-property follows by standard arguments.
Similar arguments yield analyticity of $\left(e^{tA}\right)$.
\end{proof}
\begin{Lemma}
\label{D}
Assume that $\lambda>0$ and $r>0$. Then
\[\psi_{\lambda} \in\mathrm{dom}((r - A)^{\alpha})\quad \mathrm{for\,\, all}\quad \alpha\in \left [0, \frac{1}{2} + \frac{\theta}{4} \right ).\]
In particular $D_{\lambda} \in \mathcal{L}(\mathbb{R}; \mathrm{dom}((\lambda - A)^{\alpha}))$ for all $\alpha\in \left [0, \frac{1}{2} + \frac{\theta}{4} \right )$.
\end{Lemma}
\begin{proof}
We consider the case of $\rho(\xi)=\xi^{1+\theta}$ only. The other case may be proved by similar if somewhat simpler
arguments.
\par\noindent
Note first that if $\psi_{\lambda}\in (\mathrm{dom}( A), \mathcal{H})_{2,\sigma}$ then  $\psi_{\lambda} \in \mathrm{dom}\left((r -
A)^\alpha\right)$ for all $\alpha \in (0,1-\sigma)$\footnote{$(\mathrm{dom}(A), \mathcal H)_{2,\sigma}$ denotes the real interpolation space}, see for
example Theorem 11.5.1 in \cite{MartinezAlix01}. Hence the claim will follow if we show that
$\psi_{\lambda}\in (\mathrm{dom}(A),\mathcal{H})_{2,\sigma}$ for
\begin{equation}
\label{eq:rangeofsigma}
\frac{1}{2}-\frac{\theta}{4}<\sigma<\frac{1}{2} .
\end{equation}
By Theorem 10.1 of \cite{LionsMagenes72}) $\psi_{\lambda}\in \left(\mathrm{dom}(A), \mathcal{H}\right)_{2,\sigma}$ if and only if
\begin{equation}\label{e0}
\int_0^\infty t^{2\sigma-3}\left|\left(e^{tA}-I\right)\psi_\lambda\right|^2_{\mathcal H}dt<\infty
\end{equation}
and taking into account \eqref{eq:rangeofsigma} it is enough to show that
\begin{equation}
\label{eq:int-interpol} I:=\int_0^{1} t^{2\sigma -3} \left | (e^{tA} - I) \psi_{\lambda} \right
|^2_{\mathcal{H}}<\infty .
\end{equation}
To show (\ref{eq:int-interpol}) we will use (\ref{eq:solution-convolution}) and (\ref{eq:kernel}) and the definition of
$\psi_\lambda$. Denoting by $\mathbf{N}$ the cumulative distribution function of the
standard normal distribution, we obtain
\[\begin{aligned}
I&=\int_0^1t^{2\sigma-3}\int_0^\infty \xi^{1+\theta}\left|\left(e^{tA}-I\right)\psi_\lambda(\xi)\right|^2d\xi dt\\
&=\int_0^1t^{2\sigma-3}\int_0^\infty \xi^{1+\theta}\left(\int_0^\infty \frac{e^{-\frac{(\xi-\eta)^2}{4t}}}{\sqrt{4pt}}e^{-\lambda\eta}d\eta -
\int_0^\infty \frac{e^{-\frac{(\xi+\eta)^2}{4t}}}{\sqrt{4pt}}e^{-\lambda\eta}d\eta-e^{\lambda\xi}\right)^2d\xi dt\\
&=\int_0^1t^{2\sigma-3}\int_0^\infty \xi^{1+\theta}\left(e^{-\lambda\xi}e^{\lambda^2t}\mathbf N\left(\frac{\xi}{\sqrt{2t}}-\lambda\sqrt{t}\right)
-e^{\lambda\xi}e^{\lambda^2t}\left(1-\mathbf N\left(\frac{\xi}{\sqrt{2t}}+\lambda\sqrt{t}\right)\right)-e^{\lambda\xi}\right)^2d\xi dt\\
&\le 2\left(I_1+I_2+I_3\right)
\end{aligned}\]
where $I_1$, $I_2$ and $I_3$ are respectively
\[
I_1 := \int_0^{1} t^{2\sigma -3} \int_0^{+\infty} \xi^{1+\theta} \Bigg [ e^{-\lambda \xi} \left ( e^{\lambda^2 t} -1 \right ) \mathbf{N} \left ( \frac{\xi}{\sqrt{2t}} - \lambda \sqrt{2t} \right ) \Bigg ]^2 \ud \xi \ud t
\]
\[
I_2 := \int_0^{1} t^{2\sigma -3} \int_0^{+\infty} \xi^{1+\theta} \Bigg [ e^{-\lambda \xi} \left ( \mathbf{N} \left ( \frac{\xi}{\sqrt{2t}} - \lambda \sqrt{2t} \right ) -1 \right ) \Bigg ]^2 \ud \xi \ud t
\]
\[
I_3 := \int_0^{1} t^{2\sigma -3} \int_0^{+\infty} \xi^{1+\theta} \Bigg [ e^{\lambda \xi} e^{\lambda^2 t} \left ( 1-  \mathbf{N} \left ( \frac{\xi}{\sqrt{2t}} + \lambda \sqrt{2t} \right ) \right ) \Bigg ]^2 \ud \xi \ud t
\]
Since for $t\in [0,1]$ we have $\left | e^{\lambda^2 t} - 1 \right | \leq \left ( e^{\lambda^2} - 1 \right ) t$ we find
that $I_1$ converges for every $\sigma>0$. $I_3$ can be estimated, using that the standard estimate
\[(1-\mathbf{N}(s))\leq \frac{1}{s} \frac{e^{-s^2/2}}{\sqrt{2\pi}}\]
as follows:
\[\begin{aligned}
\label{eq:stimaI3}
I_3 &\leq \int_0^{1} t^{2\sigma -3} \int_0^{+\infty} \xi^{1+\theta}  e^{2\lambda \xi} e^{2\lambda^2 t} \frac{2t}{\xi^2} e^{-\xi^2/(2t)} \ud \xi \ud t \\
&\leq C_1 \int_0^{1} t^{2\sigma -2} \int_0^{+\infty} \xi^{-1+\theta}\,  e^{2\lambda\xi -\xi^2/(2t)}\, \ud \xi \ud t\\
&= C_1 \int_0^{1} t^{2\sigma -2} \int_0^{+\infty} y^{-1+\theta}\, t^{\frac{\theta-1}{2}}\,  e^{2\lambda y \sqrt{t} -y^2/2}\, t^{1/2} \ud \xi \ud t \\
&\leq C_1 \left ( \int_0^{1} t^{2\sigma -2 + \frac{\theta}{2}}\ud t \right ) \left (  \int_0^{+\infty} y^{-1+\theta}\,   e^{2\lambda y -y^2/2}  \ud \xi \right ) <\infty
\end{aligned}\]
where the finiteness of the first term follows from (\ref{eq:rangeofsigma}). The estimate for $I_2$ can be obtained in
a similar way.
\end{proof}
\subsection{Properties of the solution of the state equation}
Let $W$ be a real Brownian motion on a probability space
$(\Omega, \mathcal{F}, \mathbb{P})$  and let $\left(\mathcal{F}_t\right)$ denote the natural filtration of $W$. We need to give a rigorous meaning to equation (\ref{eq:pde-state-equation}). To this end we will assume in the sequel that
\[\lambda_0>0\quad\mathrm{and}\quad \alpha\in \left (\frac{1}{2},\frac{1}{2}+\frac{\theta}{4} \right ),\]
are fixed.
We will denote by $D$ the operator
$D_{\lambda_0} \in \mathcal{L}(\mathbb{R}; D((\lambda_0 - A)^{\alpha}))$ and $\psi_{\lambda_0}=D_{\lambda_0}(1)$. By Proposition \ref{pr:propertiessemigroup} the semigroup $\left(e^{tA}\right)$ is analytic and therefore for any $\gamma\ge 0$
\begin{equation}
\label{eq:estimateAalphaS} \| (\lambda_0 - A)^{\gamma} e^{tA} \|_{\mathcal H}\leq M_\gamma t^{-\gamma} \qquad \text{for all } t\in
(0,T],
\end{equation}
see for example \cite{Pazy83} (Theorem 6.13 page 75). By Lemma \ref{D} the operator $B=(\lambda_0 - A)D:\mathbb R\to\mathcal H^{\alpha-1}$ is bounded\footnote{For $\beta>0$
the space $\mathcal H^{-\beta}$ is defined as a completion of $\mathcal H$ with respect to the norm
$|x|_{-\beta}=\left|\left(\lambda_0-A\right)^{-\beta}x\right|$}. Moreover, for $t>0$ the operator
\[Ae^{tA}D_{\lambda_0}=\left(\lambda_0-A\right)^{1-\alpha}e^{tA}\left(\lambda_0-A\right)^{\alpha}D_{\lambda_0}:\mathbb
R\to\mathcal H\] is bounded as well.
 We will write $e^{tA}B=Ae^{tA}D_{\lambda_0}$
Now, we reformulate equation equation (\ref{eq:pde-state-equation}), still formally, as a stochastic evolution equation in $\mathcal H$:
\begin{equation}
\label{eq:state-equation-heat}
\left\{
\begin{array}{l}
\ud x(t) = \left (  A x(t) + Bu(t) \right ) \ud t + B\ud W(t)\\
x(\tau)=x_0\in\mathcal{H}
\end{array}
\right .
\end{equation}
where the control $u$ is chosen in the set $M^2_W(\tau,T;\mathbb{R})$ of progressively measurable processes endowed with the norm
\[\|u\|_{M_W^2}^2=\mathbb{E} \int_\tau^T |u(t)|^2 \ud t <\infty.\]
 The next two results show that we can give a meaning
to (\ref{eq:state-equation-heat}).
\begin{Theorem}
\label{lm:stochconvol-heat} For all $\gamma<2\alpha-1$ the following holds.\\
(i) The operator $t\to e^{tA}B:\mathbb R\to\mathcal H$ is bounded for each $t>0$ and the function
\[t\to e^{tA}Ba\]
is continuous for every $a\in\mathbb R$.\\
(ii)
\begin{equation}
\label{eq:L2Bfinite-heat} \int_0^T s^{-\gamma} \left\|\left(\lambda_0-A\right)e^{sA}\psi_{\lambda_0}\right\|_{\mathcal H}^2 \ud s < \infty .
\end{equation}
(iii) For every $T>\tau\ge 0$ the process
\[
W_A(t)= \int_\tau^t e^{(t-s)A} B \ud W(s),\quad t\in [\tau ,T]
\]
is well defined, belongs to $C([\tau,T];L^2(\Omega;\mathcal{H}))$ and has continuous trajectories in $\mathcal H$.
\end{Theorem}
\begin{proof}
(i) It follows immediately from the definition of $B$ and Lemma \ref{D} since
\begin{equation}
\label{D0} e^{tA}Ba=a(\lambda_0 - A)^{1-\alpha} e^{tA} (\lambda_0 - A)^{\alpha}\psi_{\lambda_0}, \qquad a\in\mathbb{R}
\end{equation}
(ii) (By \ref{D0}) and (\ref{eq:estimateAalphaS}) we have for
$\alpha\in\left(\frac{1}{2},\frac{1}{2}+\frac{\theta}{4}\right)$
\[\left\|\left(\lambda_0-A\right)e^{sA}\psi_{\lambda_0}\right\|_{HS}^2
=\left|\left(\lambda_0-A\right)^{1-\alpha}e^{sA}\left(\lambda_0-A\right)^\alpha\psi_{\lambda_0}\right|^2
\le\frac{C}{s^{2(1-\alpha)}}\left|\left(\lambda_0-A\right)^\alpha\psi_{\lambda_0}\right|^2\] and the estimate
(\ref{eq:L2Bfinite-heat}) follows immediately for a certain $\gamma<2\alpha-1$.\\
(iii) Using (\ref{eq:L2Bfinite-heat}) with $\gamma=0$ we find immediately that, for every $t\geq 0$, $W_A(t)$ is well
defined and (see for example \cite{DaPratoZabczyk92} Proposition 4.5 page 91)
\begin{equation}
\label{eq:L2Bfinite}
\mathbb{E} \left |W_A(t) \right|^2_{\mathcal H} = \int_\tau^t |e^{sA} ((\lambda -A)D)|^2_{HS} \ud s < \infty .
\end{equation}
Such an estimate gives also, through standard arguments, the mean square continuity.
The continuity follows from (\ref{eq:L2Bfinite-heat}) for $\gamma>0$ using a factorization argument as in
\cite{DaPratoZabczyk93} Theorem 2.3 page 174.
\end{proof}
\begin{Lemma}
\label{lm:controlterm}
Let $T>0$ be fixed, $\lambda>0$ and $u\in M_W^2(\tau,T;\mathbb{R})$. Then the process
\[
I(t)=\int_\tau^t e^{(t-s)A}Bu(s) \ud s, \quad t\le T,
\]
is well defined,  $I\in M^2_W(\tau,T;\mathcal{H})$, and there exists $C>0$ such that
\[\mathbb E\|I\|^2_{M_W^2}\le C\|u\|_{M_W^2}^2.\]
Moreover, $I$  is in
$C(\tau,T;L^2(\Omega,\mathcal{H}))$ and has continuous trajectories.
\end{Lemma}
\begin{proof}
The first part of the Lemma follow from (\ref{D0}) by standard arguments. The mean-square continuity and continuity of
$I$ follows from (\ref{eq:estimateAalphaS}) and H\"older inequality (since $\alpha>1/2$) in the expression
\[
I(t) = \int_\tau^t \left [ (\lambda_0 - A)^{1-\alpha} e^{(t-s)A} \right ]  \left [ (\lambda_0 - A)^{\alpha} Du(s) \right ] \ud s.
\]
and the claim follows.
\end{proof}
\begin{Definition}
\label{def:solution} Let $u\in M_W^2$. An $\mathcal{H}$-valued predictable process $x$, defined on $[0,T]$ is called a
mild solution of (\ref{eq:state-equation-heat}) if
\[
\mathbb{P} \left [ \int_\tau^T |x(s)|^2 \ud s <\infty \right] =1
\]
and
\[
x(t) = e^{(t-\tau)A} x_0 + \int_\tau^t e^{(t-s)A}Bu(s) \ud s + \int_\tau^t e^{(t-s)A} B \ud W(s)
\]
\end{Definition}
\begin{Theorem}
\label{pr:regsolution} Equation (\ref{eq:state-equation-heat}) has a
unique mild solution $x\in C(\tau,T;L^2(\Omega,\mathcal{H}))$. Moreover, $x$ has continuous trajectories $\mathbb P$-a.s. If
$u=0$ then equation (\ref{eq:state-equation-heat}) defines a Markov process in $\mathcal H$.
\end{Theorem}
\begin{proof}
The properties of the stochastic convolution term come from Lemma \ref{lm:stochconvol-heat}, those of $\int_\tau^t e^{(t-s)A} Bu(s)$ from Lemma \ref{lm:controlterm}. The Markov property can be proved with standard arguments (see for example \cite{DaPratoZabczyk92} Theorem 9.8 page 249).
\end{proof}
\subsection{The approximating equation}
Let $\mathcal{I}_n\nd (n(n-A)^{-1})^2$. We will approximate $x$ using
\begin{equation}
\label{eq:def-xn}
x_n\nd \mathcal{I}_n x.
\end{equation}
We have that
\begin{equation}
\label{eq:convx_ntox}
x_n \xrightarrow{C([\tau,T];L^2(\Omega,\mathcal{H}))} x.
\end{equation}
We use it to obtain more regularity and to guarantee the existence of a strong solution and then to be able to apply the Ito's rule (Proposition \ref{pr:fundamentaleq}). From Proposition \ref{pr:regsolution} we know that $x_n\in C([\tau,T];L^2(\Omega,\mathrm{dom}(A^2)))$. We have $B_n:=\mathcal{I}_n B\in \mathcal{L}(\mathbb{R}; \mathrm{dom}(A))$ and then $B_n u\in M^2_W(\tau,T;\mathrm{dom}(A))$. Furthermore, $x_n$ satisfies the following stochastic differential equation:
\begin{equation}
\label{eq:state-eq-approx}
\left \{
\begin{array}{l}
\ud x_n(t) = \left (  A x_n(t) + B_n u(t)  \right ) \ud t + B_n \ud W(t)\\
x_n(\tau) = \mathcal{I}_n x_0
\end{array}
\right .
\end{equation}
in strong (an then mild) sense (see \cite{DaPratoZabczyk92} Section 6.1).
So we have
\begin{equation}
\label{eq:strongform}
x_n(t) = \mathcal{I}_n x_0 + \int_\tau^t A x_n(s) \ud s + \int_\tau^t B_n u(s) \ud s + \int_\tau^t B_n \ud W(s)
\end{equation}
\section{The linear quadratic problem}
Let us recall that we work under the assumption
\[\frac{1}{2}<\alpha<\frac{1}{2}+\frac{\theta}{4}.\]
 We consider another Hilbert space $Y$, an operator $C\in \mathcal{L}(\mathcal{H}; Y)$ and a
symmetric and positive $G\in\mathcal{L}(\mathcal{H};\mathcal{H})$.
 For a fixed $T>0$ we define the set of the admissible controls as $\mathcal{U_\tau}=M^2_W(\tau,T;\mathbb{R})$. We consider the linear quadratic optimal control problem governed by equation (\ref{eq:state-equation-heat}) and quadratic cost functional (to be minimized)
\begin{equation}
\label{eq:functional}
J(\tau, x_0, u):=\mathbb{E} \left [ \int_\tau^T \left(|Cx(t)|_Y^2 + |u(t)|^2_\mathbb{R}\right) \ud t + \lla G x(T) ,x(T) \rra \right ].
\end{equation}
The value function of the problem is
\[
V(\tau,x_0) := \inf_{u \in \mathcal{U}_\tau} J(\tau, x_0, u)
\]
We consider now the ``associated'' deterministic linear quadratic problem. It is characterized by the state equation
\begin{equation}\label{DET}
%\tag{DET}
\left \{
\begin{array}{l}
\dot{x}(t) = A x(t) + Bu(t)\\
x(\tau)=x_0,
\end{array}
\right .
\end{equation}
by the set of admissible controls $\mathcal{U}_{DET} := L^2(\tau,T;\mathbb{R})$ and by the functional
\[
J_{DET}(\tau,x_0, u) := \int_\tau^T\left( |Cx(t)|_Y^2 + |u(t)|^2_\mathbb{R}\right) \ud t + \lla G x(T) ,x(T) \rra.
\]
In what follows we we will use the following notations.
\begin{Notation}
\begin{itemize}
\item[] $\Sigma(\mathcal{H}) = \left \{ T\in\mathcal{L}(\mathcal{H};\mathcal{H}) \; : \; T \; \text{hermitian} \right \}$
\item[] $\Sigma^+(\mathcal{H}) = \left \{ T\in \Sigma(\mathcal{H}) \; :  \; \lla Tx,x \rra \geq 0 \; \text{for all } x\in\mathcal{H} \right \}$
\item[] $C_s([0,T];\Sigma(\mathcal{H})) = \left \{ F\colon [0,T] \to \Sigma(\mathcal{H}) \; : \; F \; \text{strongly continuous} \right \}$ \end{itemize}
\end{Notation}
Note that (\cite{BDDM2} page 137) for every $P\in C_s([0,T];\Sigma(\mathcal{H}))$
\begin{equation}
\label{eq:uniflimP}
\sup_{t\in [0,T]} \|P(t)\| <\infty.
\end{equation}
The Riccati equation formally associated with the deterministic control problem (\ref{DET}) has the form
\begin{equation}
\label{eq:riccati}
\left \{
\begin{array}{l}
P'(t) = -A^*P(t) - P(t)A^* - C^*C + P(t)ABB^*A^*P(t)\\
P(T)= G,
\end{array}
\right .
\end{equation}
but the concept of solution to this equation requires a rigorous definition. We start with some notations.
\begin{Definition}
\label{def:Csalphabeta}
We denote by $C_{s,\alpha}([0,T];\Sigma(\mathcal{H}))$  the set of all $P\in C_s([0,T];\Sigma(\mathcal{H}))$ such that
\[
\begin{array}{rl}
(i) & P(t)x \in D((\lambda_0 - A^*)^{1-\alpha}) \qquad \forall x\in\mathcal{H}, \forall t\in [0,T)\\
(ii) & V_P(t) \nd (\lambda_0 - A^*)^{1-\alpha} P(t) \in C([0,T);\mathcal{L}(\mathcal{H}))\\
(iii) & \lim_{t\to T^-} \left ( (T-t)^{1-\alpha} V_P(t) x \right ) =0 \qquad \forall x\in\mathcal{H} \\
\end{array}
\]
Given $P\in C_{s,\alpha}([0,T];\Sigma(\mathcal{H}))$, the norm $|P|_{\alpha}$ is defined as
\[
|P|_{\alpha} \nd \sup_{t\in[0,T)} \|P(t)\| +  \sup_{t\in [0,T)} (T-t)^{(1-\alpha)} \|(\lambda_0 - A^*)^{1-\alpha} P(t) \|
\]
\end{Definition}
It can be proved (see \cite{BDDM2} page 205) that $C_{s,\alpha}([0,T];\Sigma(\mathcal{H}))$, endowed with the norm
$|\cdot|_{\alpha}$, is a Banach space. We will use the notation $E\nd(\lambda_0-A)^{\alpha}D \in
\mathcal{L}(\mathbb{R};\mathcal{H})$.
\par\noindent
Note that if $|P|_{\alpha}<\infty$ then (since $\alpha>1/2$)
\begin{equation}
\label{eq:normL2forP}
|P|_{L^2(0,T;\mathcal{L}(\mathcal{H}))} < \infty
\end{equation}
\begin{Definition}\label{ricc}
We say that $P\in C_{s,\alpha}([0,T];\Sigma(\mathcal{H}))$ is a weak solution of the Riccati equation
(\ref{eq:riccati}) if for all $x,y\in \mathrm{dom}(A)$ and all $t\in (0,T)$
\begin{equation}
\label{eq:riccatiweak}
\left \{
\begin{array}{l}
\frac{\ud }{\ud t} \lla P(t)x,y  \rra = - \lla P(t) x, Ay  \rra - \lla P(t) A x,y \rra - \lla Cx, Cy \rra+ \lla E^* V_P(t) x, E^* V_P(t) y \rra  \\
P(T) = G.
\end{array}
\right .
\end{equation}
\end{Definition}
We recall now the existence and uniqueness theorem for the (\ref{eq:riccatiweak}):
\begin{Theorem}
\label{th:flandoli-existence} (i) The Riccati equation (\ref{eq:riccatiweak}) has a unique weak solution in $P$ in
$C_{s,\alpha}([0,T];\Sigma^+(\mathcal{H}))$\\
(ii) $P\in C_{s,\alpha}([0,T];\Sigma^+(\mathcal{H}))$ is a weak solution of (\ref{eq:riccati}) if and only if
it solves the following mild equation:
\begin{multline}
\label{eq:riccatimild}
P(t) = e^{(T-t)A^*} G e^{(T-t)A} + \int_t^T e^{(s-t)A^*} C^*C e^{(s-t)A} \ud s\\
+ \int_t^T e^{(s-t)A^*} V_P^*(s) E E^* V_P(s) e^{(s-t)A} \ud s
\end{multline}
\end{Theorem}
\begin{proof}
See \cite{BDDM2} Theorem 2.1 page 207 for the proof of (i) and \cite{BDDM2} Proposition 2.1 page 206 for (ii).
\end{proof}

\subsection{Dynamic Programming}

\begin{Lemma}
\label{lm:Ptracefinite}
We have that
\[
\int_0^T  \left\langle ((\lambda_0 -A) \psi_{\lambda_0}), P(t)((\lambda_0 -A) \psi_{\lambda_0}) \right \rangle_{\mathcal{H}} \ud t <\infty
\]
\end{Lemma}
\begin{proof}
We use the fact that $P$ satisfies the mild equation (\ref{eq:riccatimild}). We have that
\begin{multline}
\int_0^T \left\langle ((\lambda_0 -A) \psi_{\lambda_0}), P(t)((\lambda_0 -A) \psi_{\lambda_0}) \right \rangle_{\mathcal{H}} \ud t \\
=\int_0^T E^* (\lambda_0 -A^*)^{1-\alpha} P(t) (\lambda_0 -A)^{1-\alpha} E (1) \ud t \\
= I_1 + I_2 + I_3 \\
\nd  \int_0^T E^* (\lambda_0 -A^*)^{1-\alpha} e^{(T-t)A^*} G e^{(T-t)A} (\lambda_0 -A)^{1-\alpha} E (1) \ud t\\
+ \int_0^T E^* (\lambda_0 -A^*)^{1-\alpha} \left ( \int_t^T e^{(s-t)A^*} C^*C e^{(s-t)A} \ud s \right ) (\lambda_0 -A)^{1-\alpha} E (1) \ud t\\
+ \int_0^T E^* (\lambda_0 -A^*)^{1-\alpha} \int_t^T e^{(s-t)A^*} V_P^*(s) E E^* V_P(s) e^{(s-t)A} \ud s (\lambda_0 -A)^{1-\alpha} E (1) \ud t.
\end{multline}
For $I_1$ we have only to check the integrability for $t\to T$ and it follows from the fact that $\alpha>1/2$ and from (\ref{eq:estimateAalphaS}): $\left \| \left ( (\lambda_0 -A^*)^{1-\alpha} e^{(T-t)A^*} \right ) \right \| \leq M_{1-\alpha} (T-t)^{1-\alpha}$.
For $I_2$ we proceed in a similar way: we can write $I_2$ as:
\[
I_2= \int_0^T \int_t^T \left | C \left ( (\lambda_0 -A)^{1-\alpha} e^{(s-t)A}  \right ) E (1) \right |^2 \ud s \ud t
\]
and we can conclude as for $I_1$, using (\ref{eq:estimateAalphaS}).
For $I_3$ we can observe that:
\[
I_3 = \int_0^T  \int_t^T \left | E^* (\lambda_0 -A^*)^{1-\alpha} e^{(s-t)A^*} V_P^*(s) E(1) \right |^2\ud s   \ud t
\]
Note that from $(ii)$ of Definition \ref{def:Csalphabeta} and from the finiteness of the norm $|P|_{\alpha}$ we know that
\begin{equation}\label{anal}
\|V_P^*(s)\| \leq \frac{C_1}{(T-s)^{1-\alpha}}
\end{equation}
and
\[
\|E^* (\lambda_0 -A^*)^{1-\alpha} e^{(s-t)A^*}\| \leq \frac{C_2}{(s-t)^{1-\alpha}}
\]
The claim follows by straightforward computations.
\end{proof}
\begin{Proposition}
\label{pr:fundamentaleq}
If $u\in M^2_W(\tau,T;\mathbb{R})$ is a control and $x$ is the related trajectory, then
\begin{multline}
\label{eq:fundamenteq}
\mathbb{E} \left [ \lla G x(T),x(T) \rra + \int_\tau^T |C x(t)|^2_{Y} + |u(t)|_{\mathbb{R}}^2 \ud t \right ]  \\
= \lla P(\tau) x(\tau) ,x(\tau) \rra  + \mathbb{E} \left [ \int_\tau^T |u(t) + E^* V_P(t)x(t) |^2_{\mathbb{R}} \right ] \\
+ \int_\tau^T \frac{1}{2} \left\langle (\lambda_0 -A) \psi_{\lambda_0}, P(t)((\lambda_0 -A) \psi_{\lambda_0})\right\rangle_{\mathcal{H}} \ud t .
\end{multline}
\end{Proposition}
\begin{proof}
We will perform the following steps: we first approximate $x$ using $x_n$ defined in (\ref{eq:def-xn}) then we compute $\int_\tau^{T_0} \frac{\ud}{\ud t} \lla P(t) x_n(t), x_n(t) \rra \ud t$ using Ito's formula and eventually we will consider to the limit $n \to\infty$ and then $T_0\to T$.
Let $\mathcal{L}(\mathrm{dom}(A);\mathcal{H})$ be the space of bounded operators from $\mathrm{dom}(A)$ endowed with the graph norm to $\mathcal H$.
Note that since $P\in C([\tau,T);\mathcal{L}(\mathcal{H};\mathcal{H}))$ it is \textit{a fortiori} and element of $C([\tau,T);\mathcal{L}(\mathrm{dom}(A);\mathcal{H}))$. Consider $T_0< T$ and the following function ($\mathrm{dom}(A)$ is endowed with the graph norm)
\[
\left \{
\begin{array}{l}
\Phi\colon [\tau,T_0]\times \mathrm{dom}(A) \to \mathbb{R}\\
\Phi\colon (t,x) \mapsto \lla P(t) x, x \rra_{\mathcal{H}} .
\end{array}
\right.
\]
Note that in the definition of $\Phi$ we use the scalar product of $\mathcal{H}$ and not of $\mathrm{dom}(A)$. $\Phi$ is twice continuously differentiable with locally bounded derivatives in $x$ on $[\tau,T_0]\times \mathrm{dom}(A)$. Moreover we have that $\partial_x \Phi (t,x) = P(t) x$  and $\partial_x^2 \Phi (t,x)(y,z) = 2 \lla P(t) y, z \rra_{\mathcal{H}}$. The first derivative in $t$ is also continuous and locally bounded on $[\tau,T_0]\times \mathrm{dom}(A)$. Invoking (\ref{eq:riccatiweak}) we have
\[
\begin{array}{l}
\frac{\ud }{\ud t} \lla P(t)x,y  \rra_{\mathcal{H}} = - \lla P(t) x, Ay  \rra_{\mathcal{H}} - \lla P(t) A x,y \rra - \lla Cx, Cy \rra_{\mathcal{H}} + \\
\qquad\qquad + \lla E^* V_P(t) x, E^* V_P(t) y \rra_{\mathcal{H}}.
\end{array}
\]
Such an expression can be discontinuous for $t=T$ only (this is the reason why we have considered a $T_0<T$). We have
already observed that $x_n$ satisfy the integral equation (\ref{eq:strongform}) also in $\mathrm{dom}(A)$ and then we can use the
Ito's rule (see \cite{DaPratoZabczyk92} page 105): we have that
\begin{multline}
\lla P({T_0}) x_n({T_0}), x_n({T_0}) \rra = \lla P(\tau) x_n(\tau), x_n(\tau) \rra - \int_\tau^{T_0} \lla C x_n(t), Cx_n(t) \rra_Y \ud t\\
- 2 \int_\tau^{T_0} \lla P(t) x_n(t) , A x_n(t) \rra \ud t
+ \int_\tau^{T_0} \lla E^* V_P(t) x_n(t), E^* V_P(t) x_n(t) \rra_{\mathbb{R}} \ud t \\
+ 2 \int_\tau^{T_0} \lla V_P(t) x_n(t), \mathcal{I}_n (\lambda_0-A)^{\alpha} Du(t) \rra \ud t + 2 \int_\tau^{T_0} \lla P(t) x_n(t) , A x_n(t) \rra \ud t \\
+ 2 \int_\tau^{T_0} \lla V_P(t) x_n(t), \mathcal{I}_n (\lambda_0-A)^{\alpha} D \ud W(t) \rra \ud t \\
+ \int_\tau^{T_0} \frac{1}{2} \left\langle ((\lambda_0 -A) \mathcal{I}_n \psi_{\lambda_0}), P(t)((\lambda_0 -A) \mathcal{I}_n \psi_{\lambda_0}) \right \rangle_{\mathcal{H}} \ud t.
\end{multline}
By simplifying the terms $\lla P(t) x_n (t), Ax_n(t) \rra$, adding and subtracting $|u(t)|_{\mathbb{R}}^2$ and $2 \int_0^{T_0} \lla u(t), E^* V_P (t) x_n(t) \rra_{\mathbb{R}}$ inside the integral and taking the expectation we find:
\begin{multline}
\mathbb{E} \left [ \lla P(T_0) x_n({T_0}),x_n({T_0}) \rra + \int_\tau^{T_0} |C x_n(t)|^2_{Y} + |u(t)|_{\mathbb{R}}^2 \ud t \right ] \\
= \lla P(\tau) x_n(\tau) ,x_n(\tau) \rra  + \mathbb{E} \left [ \int_\tau^{T_0} |u(t) + E^* V_P(t)x_n(t) |^2_{\mathbb{R}} \right ] \\
+ 2 \mathbb{E} \left [ \int_\tau^{T_0} \lla V_P(t) x_n(t), (\mathcal{I}_n-I) (\lambda_0-A)^{\alpha} Du(t) \rra \ud t \right ]\\
+ \int_\tau^{T_0} \frac{1}{2} \left\langle ((\lambda_0 -A) \mathcal{I}_n \psi_{\lambda_0}), P(t)((\lambda_0 -A) \mathcal{I}_n \psi_{\lambda_0}) \right \rangle_{\mathcal{H}} \ud t.
\end{multline}
We want now to pass with $n\to\infty$. Since by (\ref{eq:convx_ntox})
we have $x_n\xrightarrow[n\to\infty]{C([\tau,T_0];L^2(\Omega,\mathcal{H}))} x$, it is clear that
\[\lim_{n\to\infty}\mathbb E\left(\left\langle P\left(T_0\right)x_n\left(T_0\right),x_n\left(T_0\right)\right\rangle+\int_\tau^{T_0}\left|Cx_n(t)\right|_Y^2dt\right)
=\mathbb E\left(\left\langle P\left(T_0\right)x\left(T_0\right),x\left(T_0\right)\right\rangle +\int_\tau^{T_0}\left|Cx(t)\right|_Y^2dt\right)\]
and
\[\lim_{n\to\infty}\left\langle P(\tau)x_n(\tau),x_n(\tau)\right\rangle=\left\langle P(\tau)x(\tau),x(\tau)\right\rangle .\]
Since  $V_P\in
C(\left[\tau,T_0\right];\mathcal{L}(\mathcal{H}))$ and $(\lambda_0 -A)^{\alpha}D=E$ is bounded, the Dominated convergence yields
\[\lim_{n\to\infty}\mathbb E\left(\int_\tau^{T_0}\left|u(t)+E^\star V_P(t)x_n(t)\right|^2_{\mathbb R}dt+2\int_\tau^{T_0}\left\langle V_P(t)x_n(t),\left(\mathcal I_n-I\right)\left(\lambda_0-A\right)^\alpha Du(t)\right\rangle dt\right)\]
\[=\mathbb E\int_\tau^{T_0}\left|u(t)+E^\star V_P(t)x(t)\right|^2_{\mathbb R}dt .\]
Finally, using the arguments similar to those in the proof of Lemma
\ref{lm:Ptracefinite} we obtain
\[\lim_{n\to\infty}\mathbb E\int_\tau^{T_0}\frac{1}{2}\left\langle\left(\lambda_0-A\right)\mathcal I_n\psi_{\lambda_0},P(t)\left(\lambda_0-A\right)\mathcal I_n\psi_{\lambda_0} \right\rangle_{\mathcal H}dt=
\mathbb E\int_\tau^{T_0}\frac{1}{2}\left\langle\left(\lambda_0-A\right)\psi_{\lambda_0},P(t)\left(\lambda_0-A\right)\psi_{\lambda_0} \right\rangle_{\mathcal H}dt\]
and therefore, putting together the above results we obtain
\begin{multline}
\label{eq:T0fondamentale}
\mathbb{E} \left [ \lla P(T_0) x({T_0}),x({T_0}) \rra + \int_\tau^{T_0} |C x(t)|^2_{Y} + |u(t)|_{\mathbb{R}}^2 \ud t \right ] \\
= \lla P(\tau) x(\tau) ,x(\tau) \rra  + \mathbb{E} \left [ \int_\tau^{T_0} |u(t) + E^* V_P(t)x(t) |^2_{\mathbb{R}} \right ] \\
+ \int_\tau^{T_0} \frac{1}{2} \left\langle ((\lambda_0 -A) \psi_{\lambda_0}), P(t)((\lambda_0 -A) \psi_{\lambda_0}) \right \rangle_{\mathcal{H}} \ud t.
\end{multline}
Now we pass to the limit in $T_0\uparrow T$ in (\ref{eq:T0fondamentale}). To show the convergence of the left hand side of \eqref{eq:T0fondamentale} it is enough to invoke monotone convergence and to show that
\begin{equation}\label{ex1}
\lim_{T_0\to T}\mathbb E\left\langle P\left(T_0\right)x\left(T_0\right),x\left(T_0\right)\right\rangle=\mathbb E\left\langle P\left(T\right)x\left(T\right),x\left(T\right)\right\rangle .
\end{equation}
To this end note that
\begin{multline}\label{ex2}
\left\langle P(T) x(T) , x(T) \right\rangle -\left\langle P\left(T_0\right) x\left(T_0\right) , x\left(T_0\right) \right\rangle \\
= \left\langle \left(P(T)-P\left(T_0\right)\right) x(T) , x(T) \right\rangle \\
+ \left\langle P(T_0) (x(T)-x(T_0)) , x(T) \right\rangle
+ \left\langle P(T_0) x(T_0) , x(T) - x(T_0) \right\rangle.
\end{multline}
Then the strong continuity of $P$ at $T$ yields
\[\lim_{T_0\to T}\left\langle \left(P(T)-P\left(T_0\right)\right) x(T) , x(T) \right\rangle=0\]
hence
by (\ref{eq:uniflimP}) and the fact that $x\in C([\tau,T];L^2(\Omega;\mathcal{H}))$ and the Dominated Convergence we obtain
\[\lim_{T_0\to T}\mathbb E\left\langle \left(P(T)-P\left(T_0\right)\right) x(T) , x(T) \right\rangle=0.\]
Again, since $x\in C([\tau,T];L^2(\Omega;\mathcal{H}))$, we find that
\[\left|\mathbb E\left\langle P(T_0) (x(T)-x(T_0)) , x(T) \right\rangle\right|\le
\sup_{t\le T}\|P(t)\|\left(\sup_{t\le T}\mathbb E|x(t)|^2\right)^{1/2}\left(\mathbb E\left|x(T)-x\left(T_0\right)\right|^2\right)^{1/2}\]
and therefore
\[\lim_{T_0\to T}\mathbb E\left\langle P(T_0) (x(T)-x(T_0)) , x(T) \right\rangle=0 .\]
By the same arguments we obtain
\[\lim_{T_0\to T}\mathbb E\left\langle P(T_0) x(T_0) , x(T) - x(T_0) \right\rangle=0\]
and therefore  we obtain the convergence of the left hand side of (\ref{eq:T0fondamentale}).
To prove convergence of the second term in the right side of (\ref{eq:T0fondamentale}) it is enough to show that $V_Px \in M^2_W(\tau,T;\mathcal{H})$. Indeed, invoking \eqref{anal} we have
\begin{equation}
\label{eq:estimateVPx}
\begin{aligned}
\mathbb{E} \int_0^T |V_P(s) x(s)|^2 \ud s &\leq C_3  \int_0^T \| V_P(s)\|^2 \mathbb{E}|x(s)|^2 \ud s \\
&\leq C_4 |x|_{C([\tau,T];L^2(\Omega;\mathcal{H}))}^2 \int_0^T \| V_P(s)\|^2 \ud s <\infty .
\end{aligned}
\end{equation}
The convergence for the third term of the right side of (\ref{eq:T0fondamentale}) for $T_0\to T$ follows from Lemma \ref{lm:Ptracefinite}.
\end{proof}
\begin{Theorem}
\label{th:optimalfeedback}
Let $\tau\in [0,T]$ and $x_0$ be in $\mathcal{H}$. Then there exists a unique optimal pair $(u^*, x^*)$ at $(\tau,x_0)$. The optimal control $u^*$ is given by the feedback formula
\begin{equation}\label{efeed}
u^*(t)= - E^* V_P(t) x^*(t)
\end{equation}
and the value function of the problem is
\[
V(\tau,x_0)= \lla P(\tau) x_0 ,x_0 \rra + \int_\tau^T \frac{1}{2} \left\langle ((\lambda_0 -A) \psi_{\lambda_0}), P(s)((\lambda_0 -A) \psi_{\lambda_0}) \right \rangle_{\mathcal{H}} \ud s
\]
\end{Theorem}
\begin{proof}
We begin proving that the equation
\begin{equation}\label{eq:mild-optimal}
x^*(t) = e^{(t-\tau)A} x_0 - \int_\tau^t (\lambda_0-A)^{1-\alpha} e^{(t-s)A} E E^* V_P(s) x^*(s) \ud s+\int_\tau^t e^{(t-s)A} B \ud W(s)
\end{equation}
has a unique solution and it is in $C([\tau,T];L^2(\Omega;\mathcal{H}))$. Consider the mapping
\[
\left \{
\begin{array}{ll}
\multicolumn{2}{l}{\phi\mapsto \Psi(\phi)}\\
\Psi(\phi) = \!\!\!& e^{(t-\tau)A} x_0\! - \displaystyle \int_\tau^t (\lambda_0-A)^{1-\alpha} e^{(t-s)A} E E^* V_P(s) \phi(s) \ud s \\
& + \displaystyle \int_\tau^t e^{(t-s)A} B \ud W(s).
\end{array}
\right .
\]
We want to prove that $\Psi(\phi)$ defines a contraction on $C([\tau,t];L^2(\Omega;\mathcal{H}))$ if we choose $t$ small enough. Consider $\psi$ and $\phi$ in $C([\tau,T];L^2(\Omega;\mathcal{H}))$:
\begin{multline}
\mathbb{E} \left [ | (\Psi(\psi)- \Psi(\phi))(t) |^2 \right ] \\ = \mathbb{E} \left [ \left |  \int_\tau^t \left ( (\lambda_0-A)^{1-\alpha} e^{(t-s)A} \right ) E E^* V_P(s) (\psi-\phi)(s) \ud s \right |^2 \right ] \\
\le C_1 \mathbb{E} \left [ |V_P|_{L^2(0,T;\mathcal{L}(\mathcal{H}))}^2 \int_\tau^t \frac{1}{(t-s)^{2(1-\alpha)}} |(\psi-\phi)(s)|^2 \ud s \right ] \\
\le C_2 |(\psi-\phi)(s)|^2_{C([\tau,t];L^2(\Omega;\mathcal{H}))}  \int_\tau^t \frac{1}{(t-s)^{2(1-\alpha)}}  \ud s
\end{multline}
where the constants $C_1$ and $C_2$ do not depend on $t$. So if $t$ is small enough $\Psi$ is a contraction on ${C([\tau,t];L^2(\Omega;\mathcal{H}))}$. Similar estimates (together with the fact that $W_A\in C(\tau,T;L^2(\Omega;\mathcal{H}))$ prove that the image of $\Psi$ is in ${C([\tau,t];L^2(\Omega;\mathcal{H}))}$. Proceeding by iterations (we can choose an uniform step) we have the existence and uniqueness of the solution of the (\ref{eq:mild-optimal}) on ${C([\tau,T];L^2(\Omega;\mathcal{H}))}$.

We will prove now that $u^\star$ defined by \eqref{efeed} is the optimal control. Its admissibility (that is  $u^* \in M_W^2(\tau,T;\mathbb{R})$) can be proved using the same argument we used in (\ref{eq:estimateVPx}).

Now we observe that Proposition \ref{pr:fundamentaleq} implies, for every $u\in M^2_W(\tau,T;\mathbb{R})$,
\begin{equation}\label{x4}
J(\tau,x_0, u) \geq \lla P(\tau) x_0 ,x_0 \rra
+ \int_\tau^T \frac{1}{2} \left\langle ((\lambda_0 -A) \psi_{\lambda_0}), P(t)((\lambda_0 -A) \psi_{\lambda_0}) \right \rangle_{\mathcal{H}} \ud t
\end{equation}
and the couple $(u^*, x^*)$ satisfies
\begin{equation}\label{x5}
J(\tau,x_0, u^*) = \lla P(\tau) x_0 ,x_0 \rra
+ \int_\tau^T \frac{1}{2} \left\langle ((\lambda_0 -A) \psi_{\lambda_0}), P(t)((\lambda_0 -A) \psi_{\lambda_0}) \right \rangle_{\mathcal{H}} \ud t
\end{equation}
so it is optimal.
If $(\hat u , \hat x)$ is another optimal couple then by \eqref{x4} and \eqref{x5} we have
\[J(\tau,x_0, \hat{u}) = \lla P(\tau) x_0 ,x_0 \rra
+ \int_\tau^T \frac{1}{2} \left\langle ((\lambda_0 -A) \psi_{\lambda_0}), P(t)((\lambda_0 -A) \psi_{\lambda_0}) \right \rangle_{\mathcal{H}} \ud t\]
and then (\ref{eq:T0fondamentale}) yields
\[\left| \hat{u}(t) + E^* V_P(t) \hat{x}(t)\right|=0\quad dt\otimes\mathbb P-a.e. \]
 and then $\hat x$ satisfies (\ref{eq:mild-optimal}) but the solution to (\ref{eq:mild-optimal}) is unique by Theorem \ref{pr:regsolution} solution and finally we can choose continuous versions of $\hat{x}$ and $\hat{u}$ such that $x^*=\hat x$ and $u^*=\hat u$.
\end{proof}

\bibliography{biblio}
\bibliographystyle{plain}

\end{document}